\documentclass[12pt]{article}
\usepackage{}
\usepackage{amssymb}
\usepackage{mathrsfs}
\usepackage{amsmath,amssymb,theorem}
\usepackage{graphicx}
\usepackage{subfigure}
\usepackage{psfrag}
\usepackage{color}
\usepackage{enumerate}
\newtheorem{thm}{Theorem}[section]

\theorembodyfont{\rmfamily}
\def\less{\backslash}

\def\qed{ \hfill $\square$}

\title{Exact rainbow numbers for matchings in plane triangulations}

\author{\small {Zhongmei Qin$^{1,2}$, Yongxin Lan$^3$, Yongtang Shi$^3$, Jun Yue$^4$}\\[5mm]
{\small $^1$School of Science}\\
{\small Chang'an University, Xi'an, Shaanxi 710064, P.R. China}\\
{\small $^2$Department of Applied Mathematics}\\
{\small Northwestern Polytechnical University, Xi'an, Shaanxi 710072, P.R. China}\\
{\small  $^3$Center for Combinatorics and LPMC}\\
{\small Nankai University, Tianjin 300071, P.R. China}\\
{\small  $^4$School of Mathematics and Statistics}\\
{\small Shandong Normal University, Jinan, Shandong 250358, P.R. China}\\
{\small Emails: qinzhongmei90@163.com, lan@mail.nankai.edu.cn,}\\
 {\small shi@nankai.edu.cn, yuejun06@126.com}\\
}
\date{}
\begin{document}\maketitle
\begin{abstract}
Given two graphs $G$ and $H$, the {\it rainbow number} $rb(G,H)$ for $H$ with respect to $G$ is defined as the minimum number $k$ such that any $k$-edge-coloring of $G$ contains a rainbow $H$, i.e., a copy of $H$, all of its edges have different colors. Denote by $M_t$ a matching of size $t$ and $\mathcal {T}_n$ the class of all plane triangulations of order $n$, respectively. Jendrol$'$,  Schiermeyer and Tu initiated to investigate the rainbow numbers for matchings in plane triangulations, and proved some bounds for the value of $rb({\mathcal {T}_n},M_t)$. Chen, Lan and Song proved that $2n+3t-14 \le rb(\mathcal {T}_n, M_t)\le 2n+4t-13$ for all $n\ge 3t-6$ and $t \ge 6$. In this paper, we determine the exact values of  $rb({\mathcal {T}_n},M_t)$ for large $n$, namely,  $rb({\mathcal {T}_n},M_t)=2n+3t-14$ for all $n \ge 9t+3$ and $t\ge 7$.\\[2mm]
\end{abstract}

\section{Introduction}

All graphs in this paper are undirected, finite and simple. We follow \cite{BM08} for graph theoretical notation and terminology not defined here. Let $G$ be a connected graph with vertex set $V(G)$ and edge set $E(G)$. For any two disjoint subsets $X$ and $Y$ of $V(G)$, we use $E_G(X,Y)$ to denote the set of edges of $G$ that have one end in $X$ and the other in $Y$. %Denote by $|E_G(X,Y)|$ the number of edges in $E_G(X,Y)$.
We also denote $E_G(X,X)=E_G(X)$.
Let $e(G)$ denote the number of edges of $G$, $e_G(X,Y)$ the number of edges of $E_G(X,Y)$, $e_G(X)$ the number of edges of $E_G(X)$.
If $X=\{x\}$, then we write $E_G(x,Y)$ and $e_G(x,Y)$, respectively. For a vertex $x\in V(G)$, we use $N_G(x)$ to denote the set of vertices in $G$ which are adjacent to $x$. We define $d_G(x)=|N_G(x)|$. Given $X, Y \subseteq V(G)$,
the subgraph of $G$ induced by $X$, denoted $G[X]$, is the graph with vertex set $X$ and edge set $\{xy \in E(G) : x, y \in X\}$. We denote by $Y \less X$ the set $Y - X$. Since every planar bipartite graph on $n\ge 3$ vertices has at most $2n-4$ edges, we will frequently use the fact that $e_G(X,Y)\le 2(|X|+|Y|)-4$ when $G$ is planar and $|X\cup Y|\ge 3$. For any  positive integer $k$, let $[k]:= \{1,2,\ldots,k\}$. We use the convention that ``$A:=$" means that $A$ is denote to be the right-hand side of the relation.

A subgraph of an edge-colored graph is {\it rainbow} if all of its edges are colored distinct. Given two graphs $G$ and $H$, the {\it rainbow number} $rb(G,H)$ for $H$ with respect to $G$ is defined as the minimum number $k$ such that any $k$-edge-coloring of $G$ contains a rainbow copy of $H$. Actually, the rainbow number is closely related to anti-Ramsey number, which was introduced by Erd\H{o}s, Simonovits and S\'{o}s \cite{ESS75} in 1975. The {\it anti-Ramsey number}, denoted by $f(K_n,H)$, is the maximum number $c$ for which there is a way to color the edges of $K_n$ with $c$ colors such that every subgraph $H$ of $K_n$ has at least two edges of the same color. Clearly, $rb(K_n, H)=f(K_n, H)+1$. The rainbow numbers for some other special graph classes in complete graphs have been investigated, see \cite{A83,ESS75,HY12,J02,JW03,JL09,MN02}. Meanwhile, the researchers studied the rainbow number or the anti-Ramsey number when host graph changed from the complete graph to others, such as complete bipartite graphs (\cite{AJK04,LTJ09}), planar graphs (\cite{HJSS15,JST14,JYNC,LSS172}), hypergraphs (\cite{OY13, GLS}), and so on. For more results on rainbow numbers or anti-Ramsey numbers, we refer the readers to the survey \cite{FMO10}.

It seems that the rainbow numbers for matchings are of special interest. Rainbow numbers for matchings in complete graphs has been completely determined step by step by Schiermeyer \cite{S04}, Chen,
Li and Tu \cite{CLT09}, and Fujita, Kaneko, Schiermeyer and Suzuki \cite{FKSS09}. Rainbow numbers for matchings in complete bipartite graphs were studied by Li, Tu and Jin \cite{LTJ09}; and in hypergraphs were studied by \"Ozkahya and Young \cite{OY13}, and Frankl and Kupavskii \cite{FK} recently.

In this paper, we study the rainbow numbers in plane triangulations, which was initiated by Hor\v{n}\'{a}k, Jendrol$'$, Schiermeyer and Sot\'{a}k \cite{HJSS15}.
Let $\mathcal{T}_n$ denote the class of all plane triangulations of order $n$. We denote by $rb(\mathcal{T}_n,H)$ the minimum number of colors $k$ such that, if $H\subseteq T_n\in \mathcal{T}_n$, then any edge-coloring of $T_n$ with at least $k$ colors contains a rainbow copy of $H$.  Hor\v{n}\'{a}k et al. \cite{HJSS15} investigated the rainbow numbers for cycles. Very recently, Lan, Shi and Song \cite{LSS172} improved some bounds for the rainbow number of cycles, and also got some results for paths.  Jendrol$'$, Schiermeyer and Tu \cite{JST14} investigated the rainbow numbers for matchings in plane triangulations. For all $t\ge 1$, let $M_t$ denote a matching of size $t$. In \cite{JST14}, the exact values of $rb(\mathcal {T}_n,M_t)$ when $t \le 4$ were determined, and lower and upper bounds for $rb(\mathcal {T}_n,M_t)$ were also established for all $t \ge 5$ and $n\ge 2t$. In \cite{QLS19}, the exact value of $rb(\mathcal {T}_n,M_5)$ was determined and an improved upper bound for $rb(\mathcal {T}_n,M_t)$ was obtained. Recently, the exact value of $rb(\mathcal {T}_n,M_6)$ and an improved lower and upper bounds for $rb(\mathcal {T}_n,M_t)$ were also obtained by Chen, Lan and Song in \cite{CLS18}. We summarize the known
results in \cite{CLS18,JST14,QLS19} as follows.

\begin{thm}[\cite{JST14}]\label{th1}
Let $n$ and $t$ be positive integers. Then
\begin{description}
\item[(1)] $rb(\mathcal {T}_n, M_2)=\begin{cases}
                         4 & n=4;\\
                         2 & n\ge 5.
                         \end{cases}$
\item[(2)] $rb(\mathcal {T}_n, M_3)=\begin{cases}
                         8 & n=6;\\
                         n+1 & n\ge 7.
                         \end{cases}$
\item[(3)] for all $n\ge 8$, $rb(\mathcal {T}_n, M_4)=2n-1$.

\item[(4)] for all $t \ge 5$ and $n\ge 2t$, $2n+2t-9 \le rb(\mathcal {T}_n, M_t)\le 2n+2t-7+2\binom{2t-2}{3}$.
\end{description}
\end{thm}

\begin{thm}[\cite{QLS19}]\label{th2}
Let $n$ and $t$ be positive integers. Then
\begin{description}
\item[(1)] for all $n\ge 11$, $ rb(\mathcal {T}_n, M_5)=2n+1$.

\item[(2)] for all $t \ge 5$ and $n\ge 2t$, $rb(\mathcal {T}_n, M_t)\le 2n+6t-16$.
\end{description}
\end{thm}

\begin{thm}[\cite{CLS18}]\label{th3}
Let $n$ and $t$ be positive integers. Then
\begin{description}
\item[(1)]for all $n\ge 30$, $rb(\mathcal {T}_n, M_6)=2n+4$.

\item[(2)] for all $t \ge 6$ and $n\ge 3t-6$, $2n+3t-14 \le rb(\mathcal {T}_n, M_t)\le 2n+4t-13$.

\end{description}
\end{thm}

In this paper, we prove the lower bound obtained in Theorem \ref{th3} is tight when n is large.
\begin{thm}\label{them2}
\begin{description}
\item[(i)] For all $n \ge 3t-6$ and $t\ge 7$, $rb(\mathcal {T}_n,M_t)\le 2n+3t-13$.

\item[(ii)] For all $n \ge 9t+3$ and $t\ge 7$,
$ rb(\mathcal {T}_n,M_t)=2n+3t-14$.
\end{description}
\end{thm}

In order to prove our main theorem, we use the concept of planar Tur\'{a}n number of graphs which firstly studied by Dowden \cite{D16}. The {\it planar Tur\'{a}n number} of $H$, denoted $ex_{\mathcal {P}}(n,H)$, is the maximum number of edges of a planar graph on $n$ vertices containing no subgraph isomorphic to $H$. Analogous to the relation between rainbow numbers and
Tur\'{a}n numbers proved in \cite{ESS75}, rainbow numbers in plane triangulations are closely related to planar Tur\'{a}n numbers \cite{LSS172}. Namely, given a planar graph $H$ and a positive integer $n \ge |H|$,
\begin{align}\label{eq9}
2 + ex_{\mathcal {P}}(n,\mathcal {H}) \le rb(\mathcal {T}_n,H) \le ex_{\mathcal {P}}(n,H)+1,
\end{align}
where $\mathcal {H}= \{H -e : e \in E(H)\}$.

We first determine the planar Tur\'{a}n numbers for matchings. Obviously, $ex_\mathcal {P}(n, M_t)=3n-6$ for $3\le n\le 2t-1$ and $t\ge 2$.

\begin{thm}\label{them1}
For all $n \ge 2t$ and $t\ge 4$,
$ex_\mathcal {P}(n, M_t)=\min\{3n-6,2n+3t-13\}$.
\end{thm}

%In addition, we improve the new upper bounds for rainbow numbers of matchings in Theorem \ref{th2}(2) and Theorem \ref{th3}(2). Clearly, we can see that $2n+3t-14= ex_{\mathcal {P}}(n,M_{t-1})+2 \le rb(\mathcal {T}_n,M_t)\le ex_{\mathcal {P}}(n,M_{t})+1=2n+3t-12$ for all $n\ge 3t-6$ and $t\ge 5$. Moreover, we prove the following result.

The following theorem will be also used in our proof. A graph $G$ is called {\it factor-critical} if $G-v$ contains a perfect matching for each $v \in V(G)$.

\begin{thm}[\cite{LP86}]\label{th4}
Given a graph $G=(V,E)$ and $|V|=n$, let $d$ be the size of a maximum matching of $G$. Then there exists a subset $S$ with $|S|\le d$ such that
$$d=\frac{1}{2}(n-(o(G-S)-|S|)),$$
where $o(H)$ is the number of components in the graph $H$ with an odd number of vertices. Moreover, each odd component of $G-S$ is factor-critical and each even component of $G-S$ has a perfect matching.
\end{thm}

\section{Proof of Theorem \ref{them1}.}
Let $T$ be a planar graph consisting of two parts $U$ and $V$ such that: (1) $T[U]$ is a plane triangulation with $t-1$ vertices and $T[V]$ is an empty graph with $n-t+1$ vertices; (2) $T[U,V]$ is a maximal planar bipartite graph. One can see that $T$ is a planar graph without $M_t$. If $n\le3t-7$, then $e(T)=e(T[U])+e_T(U,V)=3(t-1)-6+3(n-t+1)=3n-6$. If $n \ge 3t-6$, then $e(T)=e(T[U])+e_T(U,V)=3(t-1)-6+2n-4=2n+3t-13$. Thus $ex_\mathcal {P}(n, M_t) \ge e(T)=\min\{3n-6, 2n+3t-13$\}. Since $ex_\mathcal {P}(n, M_t)\le 3n-6$, we have $ex_\mathcal {P}(n, M_t)=3n-6$ when $n\le 3t-7$.

It remains to prove that $ex_\mathcal {P}(n, M_t)\le 2n+3t-13$ when $n \ge 3t-6$. Suppose $ex_\mathcal {P}(n, M_t)\ge 2n+3t-12$ for some $t \ge 4$ and $n\ge 3t-6$. Then there exists a planar graph $G$ on $n$ vertices containing no $M_t$ as a subgraph with at least $2n+3t-12$ edges. We choose such $G$ with the minimum number of vertices $n$.
By Theorem \ref{th4}, there exists an $S \subseteq V(G)$ with $s:=|S|\le t-1$ such that $q:=o(G\backslash S)=n+s+2-2t$. Let $H_1, H_2, \ldots, H_q$ be all the odd components of $G\backslash S$. We may assume that $|H_1|\le |H_2|\le \cdots \le |H_q|$. Let $r:=\max\{i:|H_i|=1\}$. Then $n=|G|\ge |S|+(|H_1|+\cdots+|H_r|)+(|H_{r+1}|+\cdots+|H_q|)\ge s+r+3(q-r)$. It follows that $r \ge n+2s-3t+3$. Let $U=V(H_1)\cup \cdots \cup V(H_r)$ and $W:=V(G)\backslash (S \cup U)$. Then $|W|=n-s-r$.
We claim that $s\le 1$. Suppose $s\ge 2$. Then $r\ge n+2s-3t+3\ge 1$ and so $s+r\ge 3$. Thus, $e_G(S,U)\le 2(s+r)-4$. Since $r \ge n+2s-3t+3$, it follows that
\begin{align*}
2n+3t-12 & \le  e(G)= e(G[W\cup S]) +e_G(S,U)\le 3(n-r)-6+2(s+r)-4\\
&=3n-r+2s-10\le 2n+3t-13,
\end{align*}
a contradiction. Thus, we have $s\le 1$.
%We next show that $r=0$. Suppose that $r\ge 1$. Then $e_G(S,U)\le r$. By the minimality of $n$, $e(G\less U)\ge e(G)-e_G(S,U) \ge 2n+3t-12-r> \min\{3(n-r)-6, 2(n-r)+3t-13\}=ex_\mathcal{P}(n-r,M_t)$. This implies that either $G\less U$ is not a planar graph, or $G\less U$ contains a copy of $M_t$ and so $G$ contains a copy of $M_t$, a contradiction. Thus $r=0$.
We claim that $G\less S$ has only one component. Suppose $W_1,\ldots,W_\ell$ are all the components of $G\less S$ with $\ell\ge 2$. If there exists one component $W_i$ with $|W_i|\le 2$, then $e_G(S,V(W_i))+e(W_i)\le \frac{3}{2}|W_i|$. By the minimality of $n$, $e(G\less V(W_i))\ge e(G)-e_G(S,V(W_i))-e(W_i) \ge 2n+3t-12-\frac{3}{2}|W_i|> \min\{3(n-|W_i|)-6, 2(n-|W_i|)+3t-13\}=ex_\mathcal{P}(n-|W_i|,M_t)$. This implies that either $G\less W_i$ is not a planar graph, or $G\less W_i$ contains a copy of $M_t$ and so $G$ contains a copy of $M_t$, a contradiction. So we may assume that $|W_j|\ge 3$ for any $j\in[\ell]$. This follows that $r=0$ and so $n\le 3t-3-2s$ since $0=r\ge n+2s-3t+3$. Then
\begin{align*}
2n+3t-12 &= e(G)\le e(G[W_1\cup S])+\cdots+e(G[W_\ell\cup S])\\
&\le 3(|W_1|+s)-6+\cdots+3(|W_\ell|+s)-6=3n+3s(\ell-1)-6\ell\\
&\le 3n-6+(\ell-1)(3s-6)\le 3n-12+3s\le 2n+3t-13,
\end{align*}
a contradiction. Thus, $G\less S$ has a copy of $M_t$ by Theorem~\ref{th4}, a contradiction. This complete the proof. \qed

%We claim that $s+r \le 2$. Suppose that $s+r \ge 3$ and so $e_G(S,U)\le 2(r+s)-4$.
%Since $r \ge n+2s-3t+3$, it follows that
%\begin{align*}
%2n+3t-12 & =  e(G)= e(G[W\cup S]) +e_G(S,U)\le 3(n-r)-6+2(s+r)-4\\
%&=3n-r+2s-10\le 2n+3t-13,
%\end{align*}
%a contradiction. Thus, we have $s+r \le 2$. We next show that $r=0$. Suppose that $r\ge 1$. Then $s\le 1$ and so $e_G(S,U)\le 1$. By the minimality of $n$, $e(G\less U)\ge e(G)-e_G(S,U) \ge 2n+3t-13> \min\{3(n-r)-6, 2(n-r)+3t-13\}=ex_\mathcal{P}(n-r,M_t)$. This implies that either $G\less U$ is not a planar graph, or $G\less U$ contains a copy of $M_t$ and so $G$ contains a copy of $M_t$, a contradiction. Thus $r=0$. If $s+r=2$, then $s=2$ and so $n\le 3t-7$ since $0=r\ge n+2s-3t+3$, a contradiction. Now we assume $s+r\le1$. Then $n\le 3t-3-2s$ since $0=r\ge n+2s-3t+3$. We claim that $G\less S$ has only one component. Suppose $W_1,\ldots,W_\ell$ are all the components of $G\less S$ with $\ell\ge 2$. Then
%\begin{align*}
%2n+3t-12 &= e(G)\le e(G[W_1\cup S])+\cdots+e(G[W_\ell\cup S])\\
%&\le 3(|W_1|+s)-6+\cdots+3(|W_\ell|+s)-6=3n+3s(\ell-1)-6\ell\\
%&\le 3n-6+(\ell-1)(3s-6)\le 3n-12+3s\le 2n+3t-13,
%\end{align*}
%a contradiction. Thus, $G\less S$ has a copy of $M_t$ by Theorem~\ref{th4}, a contradiction. This complete the proof. \qed

\section{Proof of Theorem \ref{them2}(i).}
From Theorem \ref{them1} and Eq. \eqref{eq9}, we can get that for all $n\ge 3t-6$ and $t\ge 5$,
\begin{align*}
2n+3t-14= ex_{\mathcal {P}}(n,M_{t-1})+2 \le rb(\mathcal {T}_n,M_t)\le ex_{\mathcal {P}}(n,M_{t})+1=2n+3t-12.
\end{align*}
So now we only need to improve the upper bound. Suppose $rb(\mathcal {T}_n, M_t)\ge 2n+3t-12$. Then there exists a plane triangulation $T$ on $n$ vertices containing no rainbow $M_t$ under an edge-coloring $c$ with $2n+3t-13$ colors. Let $G$ be a rainbow spanning subgraph of $T$ with $2n+3t-13$ edges.
By Theorem \ref{them1}, $G$ contains a copy of $M_{t-1}$. Clearly, $G$ has no copy of $M_t$ because $T$ has no rainbow copy of $M_t$ under $c$.  By Theorem \ref{th4}, there exists an $S \subseteq V(G)$ with $s:=|S|\le t-1$ such that $q:=o(G\backslash S)=n+s+2-2t$. Let $H_1, H_2, \ldots, H_q$ be all the odd components of $G\backslash S$. We may assume that $|H_1|\le |H_2|\le \cdots \le |H_q|$. Let $r:=\max\{i:|H_i|=1\}$. Then $n=|G|\ge |S|+(|H_1|+\cdots+|H_r|)+(|H_{r+1}|+\cdots+|H_q|)\ge s+r+3(q-r)$. It follows that $r \ge n+2s-3t+3$. Let $V(H_i)=\{u_i\}$ for all $i \in [r]$. We may further assume that $d_G(u_1) \ge d_G(u_2) \ge \cdots \ge d_G(u_r)$. Let $U:=\{u_1,\ldots, u_r\}$ and $W:=V(G)\backslash (S \cup U)$. Then $|W|=n-s-r$.

%\textcolor{blue}{\noindent{\bf Claim:} If $s=2$ and $r=n+2s-3t+3$, then $e(G)<2n+3t-14$.}\\
%
%\textcolor{blue}{\noindent{\bf Proof:} Suppose $e(G)\ge2n+3t-14$. Note that $q-r=(n+2-2t+2)-(n-3t+7)=t-3\ge 4$ and $e_G(S,U)\le 2r$. By Theorem~\ref{th4}, we see that $H_i\cong K_3$ for any $r+1\le i\le q$ because $r=n+2s-3t+3$. We next show that $e_G(S, V(H_i))\le 5$ for any $r+1\le i\le q$. Suppose $e_G(S, V(H_j))= 6$ for some $r+1\le j\le q$. Then $e(G[S\cup V(H_j)])=e(G[S])+e_G(S,V(H_j))+e(G[V(H_j)])=e(G[S])+6+3\le 9$, which implies $e(G[S])=0$. We claim that $e_G(S, V(H_i))\le 3$ for any $r+1\le i\neq j\le q$. Suppose $e_G(S,V(H_i))\ge 4$ for some $r+1\le i\le q$ and $i\ne j$. Then $N_G(V(H_i))=S$. It follows that $G$ contains a $K_5$-minor, which contradicts the planarity of $G$. Thus, $e_G(S, V(H_i))\le 3$ for any $r+1\le i\neq j\le q$ and so $e(G)=e(G[W\cup S])+e_G(S,U)=e(G[S])+e_G(S, W)+e(G[W])+e_G(S,U)\le 0+6+3(t-4)+3(t-3)=6t-15+2r< 2n+3t-14$, a contradiction.  Then, $e_G(S, V(H_i))\le 5$ for any $r+1\le i\le q$ and so $e(G)=e(G[W\cup S])+e_G(S,U)=e(G[S])+e_G(S, W)+e(G[W])+e_G(S,U)\le 1+5(t-3)+3(t-3)=8t-23+2r< 2n+3t-14$, a contradiction.}

We first claim $s\ge3$. Suppose $s\le 2$. If $s=2$, then $r\ge n+2s-3t+3\ge 1$ as $n\ge 3t-6$. Note that $e_G(S,U)\le 2r$. Since $r\ge n+2s-3t+3=n-3t+7$, we see that
\begin{align*}
2n+3t-13 = e(G)&= e(G[W\cup S]) +e_G(S,U)\le 3(n-r)-6+2r\\
&=3n-r-6\le 2n+3t-13,
\end{align*}
which implies that $r=n-3t+7$ and $G[W\cup S]\in \mathcal{T}_{n-r}$. Moreover, $|W|=3t-9$ and $q-r=(n+2-2t+2)-(n-3t+7)=t-3\ge 4$. By Theorem~\ref{th4}, we see that $H_i\cong K_3$ for any $r+1\le i\le q$. We next show that $e_G(S, V(H_i))\le 5$ for any $r+1\le i\le q$. Suppose $e_G(S, V(H_j))= 6$ for some $r+1\le j\le q$. Then $e(G[S\cup V(H_j)])=e(G[S])+e_G(S,V(H_j))+e(G[V(H_j)])=e(G[S])+6+3\le 9$, which implies $e(G[S])=0$. We claim that $e_G(S, V(H_i))\le 3$ for any $r+1\le i\neq j\le q$. Suppose $e_G(S,V(H_i))\ge 4$ for some $r+1\le i\le q$ and $i\ne j$. Then $N_G(V(H_i))=S$. It follows that $G$ contains a $K_5$-minor, which contradicts the planarity of $G$. Thus, $e_G(S, V(H_i))\le 3$ for any $r+1\le i\neq j\le q$ and so $e(G[W\cup S])=e(G[S])+e_G(S, W)+e(G[W])\le 0+6+3(t-4)+3(t-3)=6t-15< 3(3t-7)-6$, a contradiction.  Then, $e_G(S, V(H_i))\le 5$ for any $r+1\le i\le q$ and so $e(G[W\cup S])=e(G[S])+e_G(S, W)+e(G[W])\le 1+5(t-3)+3(t-3)=8t-23< 3(3t-7)-6$, a contradiction. Hence, we may assume that $s\le1$.
%We next show that $r=0$. Suppose that $r\ge 1$. Then $e_G(S,U)\le r$ and so $e(G\less U)\ge e(G)-e_G(S,U) \ge 2n+3t-13-r>ex_\mathcal{P}(n-r,M_t)=\min\{3(n-r)-6,2(n-r)+3t-13\}$.
%This implies that either $G\less U$ is not a planar graph or $G\less U$ contains a copy of $M_t$ and so $G$ contains a copy of $M_t$. Thus $r=0$ and so $n\le 3t-3-2s$ since $0=r\ge n+2s-3t+3$.
We claim that $G\less S$ has only one component. Suppose $W_1,\ldots,W_\ell$ are all the components of $G\less S$ with $\ell\ge 2$. If there exists one component $W_i$ with $|W_i|\le 2$, then $e_G(S,V(W_i))+e(W_i)\le \frac{3}{2}|W_i|$. By the minimality of $n$, $e(G\less V(W_i))\ge e(G)-e_G(S,V(W_i))-e(W_i) \ge 2n+3t-13-\frac{3}{2}|W_i|> \min\{3(n-|W_i|)-6, 2(n-|W_i|)+3t-13\}=ex_\mathcal{P}(n-|W_i|,M_t)$. This implies that either $G\less W_i$ is not a planar graph, or $G\less W_i$ contains a copy of $M_t$ and so $G$ contains a copy of $M_t$, a contradiction. So we may assume that $|W_j|\ge 3$ for any $j\in[\ell]$. This follows that $r=0$ and so $n\le 3t-3-2s$ since $0=r\ge n+2s-3t+3$. Then $2n+3t-13 = e(G)\le e(G[W_1\cup S])+\cdots+e(G[W_\ell\cup S])\le 3(|W_1|+s)-6+\cdots+3(|W_\ell|+s)-6=3n+3s(\ell-1)-6\ell\le 3n-6+(\ell-1)(3s-6)\le 3n-12+3s\le 2n+3t-14$, a contradiction.  Then $G\less S$ has a copy of $M_t$ by Theorem~\ref{th4}, a contradiction. Thus, $s\ge 3$ and $r\ge n+2s-3t+3\ge 3$. Then, $e_G(S,U)\le 2(r+s)-4$.
%Then \textcolor{blue}{$|W|=n-s-r$} %$w:=|W|=n-s-r$
%and $e_G(S,U)\le 2(r+s)-4$ \textcolor{blue}{because} $s+r \ge 3$.
Since $r \ge n+2s-3t+3$, it follows that
\begin{align*}
2n+3t-13 = e(G)&= e(G[W\cup S]) +e_G(S,U)\le 3(n-r)-6+2(s+r)-4\\
&=3n-r+2s-10\le 2n+3t-13,
\end{align*}
which implies that $e(G)=2n+3t-13$ and $r=n+2s-3t+3$. Moreover, $G[W\cup S] \in \mathcal {T}_{n-r}$ and $e_G(S,U)=2(s+r)-4$. Then $d_G(u_1)\le 3$.
We claim that $d_G(u_r)\ge2$. Suppose $d_G(u_r)\le 1$. Then $e(G-u_r)=e(G)-d_G(u_r)\ge 2n+3t-14 > ex_{\mathcal {P}}(n-1, M_t)=2(n-1)+3t-13$, which implies that $G-u_r$ contains a $M_t$, a contradiction.
Let $j:=\max\{i\in [r]:d_G(u_i)=3\}$ and $U_1=\{u_1, u_2, \ldots, u_j\}$, $U_2=U\backslash U_1$. Then $j=2s-4\ge 2$ since $2(s+r)-4=e_G(S,U)=3j+2(r-j)$ and $s\ge 3$.
We next show that $G[S] \in \mathcal {T}_{s}$. Suppose that $G[S] \notin \mathcal {T}_{s}$. Then $e(G[S])\le3s-7$. Since $r=n+2s-3t+3$, it follows that $|H_{r+1}|=\cdots=|H_q|=3$ and $n=|G|=|S|+(|H_1|+\cdots+|H_r|)+(|H_{r+1}|+\cdots+|H_q|)=s+r+3(q-r)$.
%Then $e(G[W])\le 3(q-r)$.
By Theorem~\ref{th4}, $H_{r+1}\cong\cdots\cong H_q\cong K_3$ and so $e(G[W])=3(q-r)$.
Notice that $e_G(S,V(H_i))\le 6$ for all $r+1 \le i \le q$ because $e_G(S, U \cup V(H_i))\le 2(s+r+3)-4$ for all $r+1 \le i \le q$ and $e_G(S,U)=2(s+r)-4$.  Then $e_G(S,W)\le 6(q-r)$. Thus $e(G[S\cup W])=e(G[S])+e_G(S,W)+e(G[W])\le (3s-7)+6(q-r)+3(q-r)=3s-7+9(q-r)=3(s+3(q-r))-7=3(n-r)-7$. But $G[S\cup W]\in \mathcal {T}_{n-r}$ and so $e(G[S\cup W])=3(n-r)-6$, a contradiction.  Thus, $G[S] \in \mathcal {T}_{s}$ and so $G[S]$ has exactly $2s-4$ $3$-faces.
We claim that $W=\emptyset$. Suppose $W\neq\emptyset$ and so $q\ge r+1$. Then, $e_G(S,W)=e(G[S\cup W])-e(G[S])-e(G[W])=3(s+w)-6-(3s-6)-3(q-r)=6(q-r)$. It follows that $e_G(S,V(H_i))=6$  and so $|N_S(V(H_i))|\ge3$ for $r+1\le i\le q$. But then $G$ contains a $K_{3,3}$-minor since $|U_1|=2s-4$, a contradiction. Thus, $W=\emptyset$.
% Since $j=|U_1|=2s-4$, it follows that $W=\emptyset$.
Then $n=|G|=|S|+|U|=s+r=s+n+2s+3-3t$. Thus, $s=t-1$ and so $r=n-t+1$. Notice that $e(T)-e(G)>1$ because $T\in\mathcal{T}_{n}$ and $G$ contains no copy of $M_t$.  We also can see that $T[U_1]$ contains no edges because $G[S] \in \mathcal {T}_{s}$. It follows that there is at least one edge in $T[U_2]$ or $T[U_1,U_2]$. Assume that $u_mu_\ell\in E(T[U_2])\cup E(T[U_1,U_2])$.
Then $d_G(u_m)+d_G(u_{\ell})\le 5$ because either $u_m\in U_2$ or $u_\ell\in U_2$. Let $G'$ be the graph obtained from $G$ by deleting the vertices $u_m, u_{\ell}$ and the edge whose color is $c(u_mu_\ell)$.
Then $e(G')\ge 2n+3t-13-5-1=2n+3t-19>2n+3t-20=ex_{\mathcal {P}}(n-2, M_{t-1})$ since $n\ge 3t-6$. By Theorem \ref{them1}, there exists a rainbow $H=M_{t-1}$ in $G'$ which contains no color $c(u_mu_\ell)$, thus we obtain a rainbow $M_t=H\cup \{u_mu_{\ell}\}$ in $T$, a contradiction. The proof is complete.
\qed\\

\section{Proof of Theorem \ref{them2}(ii).}
We only need to show that $rb(\mathcal {T}_n, M_t)\le 2n+3t-14$ for $n\ge 9t+3$ and $t\ge 7$. Suppose $rb(\mathcal {T}_n, M_t)\ge 2n+3t-13$. Then there exists a plane triangulation $T$ on $n$ vertices containing no rainbow $M_t$ under an edge-coloring $c$ with $2n+3t-14$ colors. Let $G$ be a rainbow spanning subgraph of $T$ with $2n+3t-14$ edges.
By Theorem \ref{them1}, $G$ contains a copy of $M_{t-1}$. Clearly, $G$ has no copy of $M_t$ since $T$ has no rainbow copy of $M_t$ under $c$.  By Theorem \ref{th4}, there exists an $S \subseteq V(G)$ with $s:=|S|\le t-1$ such that $q:=o(G\backslash S)=n+s+2-2t$. Let $H_1, H_2, \ldots, H_q$ be all the odd components of $G\backslash S$. We may assume that $|H_1|\le |H_2|\le \cdots \le |H_q|$. Let $r:=\max\{i:|H_i|=1\}$. Then $n=|G|\ge |S|+(|H_1|+\cdots+|H_r|)+(|H_{r+1}|+\cdots+|H_q|)\ge s+r+3(q-r)$. It follows that $r \ge n+2s-3t+3\ge3$. Let $V(H_i)=\{u_i\}$ for all $i \in [r]$.
We may further assume that $d_G(u_1) \ge d_G(u_2) \ge \cdots \ge d_G(u_r)$.
Let $U:=\{u_1,\ldots, u_r\}$, $U_1:=\{u\in U| d_G(u)\ge 3\}$ and $U_2=U\backslash U_1$. Let $|U_1|=a$. Then $|U_2|=r-a$. Let $W:=V(G)\backslash (S \cup U)$. Then $w:=|W|=n-s-r$ and $e_G(S,U)\le 2(s+r)-4$ since $s+r \ge 3$. Let $B$ denote the set of vertices of all the even components of $G-S$. We first prove the following claim.\medskip

\noindent {\bf Claim 1. } For any $u_mu_\ell\in E(T)\less E(G)$, $d_G(u_m)+d_G(u_\ell)\ge 5$.\medskip

\noindent {\bf Proof.} Suppose $d_G(u_m)+d_G(u_\ell)\le 4$. Let $G'$ be the graph obtained from $G$ by deleting the vertices $u_m, u_{\ell}$ and the edge whose color is $c(u_mu_\ell)$. Then $e(G')\ge 2n+3t-14-4-1=2n+3t-19>ex_{\mathcal {P}}(n-2, M_{t-1})=2n+3t-20$. By Theorem \ref{them1}, there exists a rainbow $H=M_{t-1}$ in $G'$ which contains no color $c(u_mu_\ell)$. Thus, we obtain a rainbow $M_t=H\cup \{u_mu_{\ell}\}$ in $T$, a contradiction.\qed\medskip

We see that
\begin{align*}
2n+3t-14 &= e(G)= e(G[W\cup S]) +e_G(S,U)\\
&\le 3(n-r)-6+2(r+s)-4=3n-r+2s-10,
\end{align*}
which implies that $r\le n+2s-3t+4$. Note that $r \ge n+2s-3t+3$, we distinguish two cases to finish the proof of the theorem.\medskip

{\bf Case 1.} $r=n+2s-3t+3$.

In this case, $|H_{r+1}|=\cdots=|H_q|=3$, $|B|=0$ and $|G|=|S|+(|H_1|+\cdots+|H_r|)+(|H_{r+1}|+\cdots+|H_q|)$. Notice that $q-r=(n+s+2-2t)-(n+2s-3t+3)=t-1-s$.
By Theorem \ref{th4}, we see that $H_{r+1}\cong\cdots\cong H_q\cong K_3$.

We claim $s\ge 3$. Suppose $s\le 2$. If $s=2$, then $q-r=t-1-2\ge 4$. We next show that $e_G(S, V(H_i))\le 5$ for any $r+1\le i\le q$. Suppose $e_G(S, V(H_j))= 6$ for some $r+1\le j\le q$. Then $e(G[S\cup V(H_j)])=e(G[S])+e_G(S,V(H_j))+e(G[V(H_j)])=e(G[S])+6+3\le 9$, which implies $e(G[S])=0$. We claim that $e_G(S, V(H_i))\le 3$ for any $r+1\le i\neq j\le q$. Suppose $e_G(S,V(H_i))\ge 4$ for some $r+1\le i\le q$ and $i\ne j$. Then $N_G(V(H_i))=S$. It follows that $G$ contains a $K_5$-minor, which contradicts the planarity of $G$. Thus, $e_G(S, V(H_i))\le 3$ for any $r+1\le i\neq j\le q$ and so $e(G)=e(G[S])+e_G(S, W)+e(G[W])+e_G(S,U)\le 0+6+3(t-4)+3(t-3)+2r=2n-1< 2n+3t-14$, a contradiction.  Then, $e_G(S, V(H_i))\le 5$ for any $r+1\le i\le q$ and so $e(G)=e(G[S])+e_G(S, W)+e(G[W])+e_G(S,U)\le 1+5(t-3)+3(t-3)+2r=2n+2t-9< 2n+3t-14$, a contradiction. Hence, we may assume that $s\le1$. Then $e(G)=e(G[S])+e_G(S, U\cup W)+e(G[W])\le n+3(t-1-s)<2n+3t-14$, a contradiction.

Thus $s\ge 3$. Then, $e(G[S])\ge 3s-7$, otherwise $e(G)=e(G[S])+e_G(S,W\cup U)+e(G[W])\le 3s-8+2n-4+3(q-r)=2n+3t-15$. If $s=3$, then $a \le 2=2s-4$ since $d_G(u)\ge 3$ for any $u \in U_1$, otherwise $G$ contains a $K_{3,3}$. If $s\ge 4$, then $G[S]$ has $2s-4$ 3-faces when $e(G[S])=3s-6$ and $G[S]$ has one $4$-face and $2s-6$ $3$-faces when $e(G[S])=3s-7$. Thus, $a \le 2s-4$. Notice that $e(G[W\cup S])\ge 3(n-r)-7=3(3t-2s-3)-7$ since $2n+3t-14=e(G)=e(G[W\cup S]) +e_G(S,U) \le e(G[W\cup S])+2(s+r)-4$ and $r=n+2s-3t+3$. Then $e(G[W\cup S\cup U_1])\ge 3(3t-2s-3+a)-7$. Thus, $G[W\cup S\cup U_1]$ has $2(3t-2s-3+a)-4$ 3-faces when $e(G[W\cup S\cup U_1])= 3(3t-2s-3+a)-6$ and $G[W\cup S\cup U_1]$ has $2(3t-2s-3+a)-6$ 3-faces and one 4-face when $e(G[W\cup S\cup U_1])= 3(3t-2s-3+a)-7$. It follows that $|U_2|=r-a>2(3t-2s-3+a)-4$ since
$r=n+2s+3-3t\ge 6t+2s+6> 6t-4s+3a-10=2(3t-2s-3+a)-4+a$ when $n\ge 9t+3$ and $a\le 2s-4$. Thus, $T[U_2]$ contains at least one edge since $d_T(u)\ge 3$ for any $u \in U_2$.
By Claim 1, $T[U_2]$ contains no edge because $d_G(u_m)+d_G(u_{\ell})\le 4$ for $u_{m}u_{\ell}\in E(T[U_2])$, a contradiction.
\medskip

{\bf Case 2.} $r=n+2s-3t+4$.

In this case, $q-r=(n+s+2-2t)-(n+2s-3t+4)=t-2-s$. %Then $s\le t-4$ because $r\le q-2$.
Since $n=|G|=|S|+(|H_1|+\cdots+|H_r|)+(|H_{r+1}|+\cdots+|H_q|)+|B|=s+r+(|H_{r+1}|+\cdots+|H_q|)+|B|$, it follows that either $|H_q|=\cdots=|H_{r+1}|=3$ and $|B|=2$, or $|H_q|=5$, $|H_{q-1}|=\cdots=|H_{r+1}|=3$ and $|B|=0$.

We first assume $|H_q|=\cdots=|H_{r+1}|=3$ and $|B|=2$. If $s\ge 3$, then $e(G)=e(G[S])+e_G(S, U\cup W)+e(G[W])\le 3s-6+2n-4+3(q-r)+1=2n+3t-15$, a contradiction. If $s\le 1$, then $e(G)=e(G[S])+e_G(S, U\cup W)+e(G[W])\le n+3(t-2-s)+1<2n+3t-14$, a contradiction. Thus, $s=2$. If $e(G[S])=0$, then $e(G)=e_G(S,W)+e_G(S,U)+e(G[W])\le 6(t-4)+2(r+2)+3(t-4)+1=2n+3t-15<2n+3t-14$, a contradiction. If $e(G[S])=1$, then $e_G(S, V(H_i))\le 5$ for any $r+1\le i\le q$, otherwise $G$ contains a copy of $K_5$. Thus, $e(G)=e(G[S])+e_G(S,W)+e_G(S,U)+e(G[W])\le 1+5(t-4)+2(r+2)+3(t-4)+1=2n+2t-10<2n+3t-14$, a contradiction.

Now we assume $|H_q|=5$, $|H_{q-1}|=\cdots=|H_{r+1}|=3$ and $|B|=0$. Then $q-r-1\ge 0$. It follows that
\begin{align*}
2n+3t-14 &= e(G)= e(G[S\cup V(H_q)]) +e_G(S,U\cup (W\backslash V(H_q)))+e(G[W\backslash V(H_q)])\\
&\le 3(s+5)-6+2(n-5)-4+3(q-r-1)= 2n+3t-14,
\end{align*}
which implies that $G[S\cup V(H_q)] \in \mathcal {T}_{s+5}$ and $e_G(S,U\cup (W\backslash V(H_q)))=2(n-5)-4$. Notice that $G[S\cup V(H_q)]$ has $2(s+5)-4$ 3-faces. Then $a\le 2(s+5)-4$ since $d_G(u)\ge 3$ for any $u \in U_1$. We claim that $e_G(S, V(H_i))=6$ for any $r+1\le i\le q-1$. If $e_G(S,V(H_i))\ge 7$ for some $r+1\le i\le q-1$, then $e(G[S\cup V(H_q)\cup V(H_i)])=e(G[S\cup V(H_q)])+e_G(S,V(H_i))+e(G[V(H_i)])\ge 3(s+5)-6+7+3>3(s+8)-6$, a contradiction. Thus, $e_G(S, V(H_i))\le 6$ for any $r+1\le i\le q-1$. If $e_G(S,V(H_i))\le 5$ for some $r+1\le i\le q-1$, then $e_G(S,U\cup W\backslash V(H_q))=e_G(S,W\less V(H_q))+e_G(S,U)\le 5+6(q-r-2)+2(s+r)-4< 2(n-5)-4$, a contradiction. Thus, $e_G(S, V(H_i))=6$ for any $r+1\le i\le q-1$. Then $e(G[S\cup W\cup U_1])=3[s+5+3(q-r-1)+a]-6$. Thus, $G[S\cup W\cup U_1]$ has $2[s+5+3(q-r-1)+a]-4$ $3$-faces.  It follows that $|U_2|=r-a> 2[s+5+3(q-r-1)+a]-4$ since $r=n+2s+4-3t\ge 6t+7+2s> 6t-12-4s+3a=2[s+5+3(q-r-1)+a]-4+a$ when $n\ge 9t+3$ and $a \le 2s+6$. Thus, $T[U_2]$ contains at least one edge since $d_T(u)\ge 3$ for any $u \in U_2$. By Claim 1, $T[U_2]$ contains no edge since $d_G(u_m)+d_G(u_{\ell})\le 4$ for $u_{m}u_{\ell}\in E(T[U_2])$, a contradiction.
\qed\\

\noindent{\bf Acknowledgments.}
Zhongmei Qin was partially supported by the Fundamental Research Funds for the Central Universities (No. 300102128104). Yongxin Lan and Yongtang Shi are partially supported by National Natural Science Foundation of China,
Natural Science Foundation of Tianjin (No. 17JCQNJC00300), the China-Slovenia bilateral project ``Some topics in modern graph theory" (No.~12-6),
Open Project Foundation of Intelligent Information Processing Key Laboratory of Shanxi Province (No. CICIP2018005),
and the Fundamental Research Funds for the Central Universities.  Jun Yue was partially supported by the National Natural Science Foundation of China (No. 11626148 and 11701342)
and the Natural Science Foundation of Shandong Province (No. ZR2016AQ01).

\end{document}